\documentclass[a4paper,12pt]{amsart}
\usepackage{amssymb}
\usepackage{ifthen}
\nonstopmode \numberwithin{equation}{section}
\newtheorem{thm}{Theorem}
\newtheorem{cor}{Corollary}
\newtheorem{lem}{Lemma}

\newtheorem{conj}{Conjecture}

\theoremstyle{definition}
\newtheorem{defn}{Definition}[section]

\newtheorem{prob}[equation]{Problem}

\newenvironment{rem}{%
\bigskip
\noindent \textsl{{\sl Remark. }}}{\bigskip}
\newenvironment{rems}{%
\bigskip
\noindent \textsl{{\sl Remarks. }}}{\bigskip}

\newcounter {own}
\def\theown {\thesection       .\arabic{own}}

\newenvironment{pf}[1][]{%
 \vskip 3mm
 \noindent
 \ifthenelse{\equal{#1}{}}%
  {{\slshape Proof. }}%
  {{\slshape #1.} }%
 }%
{\qed\bigskip}

\newcounter{alphabet}
\newcounter{tmp}
\newenvironment{Thm}[1][]{\refstepcounter{alphabet}%
\bigskip%
\noindent%
{\bf Theorem \Alph{alphabet}}%
\ifthenelse{\equal{#1}{}}{}{ (#1)}%
{\bf .} \itshape}{\vskip 8pt}

\newcommand{\ID}{{\mathbb D}}

\newcommand{\IC}{{\mathbb C}}

\def\be{\begin{equation}}
\def\ee{\end{equation}}

\newcommand{\bee}{\begin{enumerate}}
\newcommand{\eee}{\end{enumerate}}

\newcommand{\blem}{\begin{lem}}
\newcommand{\elem}{\end{lem}}
\newcommand{\bthm}{\begin{thm}}
\newcommand{\ethm}{\end{thm}}
\newcommand{\bcor}{\begin{cor}}
\newcommand{\ecor}{\end{cor}}
\newcommand{\beg}{\begin{examp}}
\newcommand{\eeg}{\end{examp}}
\newcommand{\begs}{\begin{examples}}
\newcommand{\eegs}{\end{examples}}
\newcommand{\bdefe}{\begin{defn}}
\newcommand{\edefe}{\end{defn}}
\newcommand{\bprob}{\begin{prob}}
\newcommand{\eprob}{\end{prob}}
\newcommand{\bei}{\begin{itemize}}
\newcommand{\eei}{\end{itemize}}

\newcommand{\bcon}{\begin{conj}}
\newcommand{\econ}{\end{conj}}
\newcommand{\bcons}{\begin{conjs}}
\newcommand{\econs}{\end{conjs}}
\newcommand{\bprop}{\begin{propo}}
\newcommand{\eprop}{\end{propo}}
\newcommand{\br}{\begin{rem}}
\newcommand{\er}{\end{rem}}
\newcommand{\brs}{\begin{rems}}
\newcommand{\ers}{\end{rems}}
\newcommand{\bo}{\begin{obser}}
\newcommand{\eo}{\end{obser}}
\newcommand{\bos}{\begin{obsers}}
\newcommand{\eos}{\end{obsers}}
\newcommand{\bpf}{\begin{pf}}
\newcommand{\epf}{\end{pf}}
\newcommand{\ba}{\begin{array}}
\newcommand{\ea}{\end{array}}
\newcommand{\beq}{\begin{eqnarray}}
\newcommand{\beqq}{\begin{eqnarray*}}
\newcommand{\eeq}{\end{eqnarray}}
\newcommand{\eeqq}{\end{eqnarray*}}

\newcommand{\ra}{\rightarrow}

\newcounter{minutes}
\divide\time by 60
\newcounter{hours}
\multiply\time by 60 \addtocounter{minutes}{-\time}
\begin{document}
\title{On the Taylor coefficients of a subclass of meromorphic univalent functions }
\begin{center}
{\tiny \texttt{FILE:~\jobname .tex,
        printed: \number\year-\number\month-\number\day,
        \thehours.\ifnum\theminutes<10{0}\fi\theminutes}
}
\end{center}

\author{Bappaditya Bhowmik${}^{~\mathbf{*}}$}
\address{Bappaditya Bhowmik, Department of Mathematics,
Indian Institute of Technology Kharagpur, Kharagpur - 721302, India.}
\email{bappaditya@maths.iitkgp.ernet.in}
\author{Firdoshi Parveen}
\address{Firdoshi Parveen, Department of Mathematics,
Indian Institute of Technology Kharagpur, Kharagpur - 721302, India.}
\email{frd.par@maths.iitkgp.ernet.in}

\begin{abstract}
 Let $\mathcal{V}_p(\lambda)$ be the collection of all functions $f$ defined in the unit disc $\ID$ having a
 simple pole at $z=p$ where $0<p<1$ and analytic in $\ID\setminus\{p\}$ with $f(0)=0=f'(0)-1$ and satisfying the differential inequality
 $|(z/f(z))^2 f'(z)-1|< \lambda $ for $z\in \ID$, $0<\lambda\leq 1$. Each $f\in\mathcal{V}_p(\lambda)$ has the following Taylor expansion:
 $$
 f(z)=z+\sum_{n=2}^{\infty}a_n(f) z^n, \quad |z|<p.
 $$
 In \cite{BF-3}, we conjectured that
 $$
 |a_n(f)|\leq \frac{1-(\lambda p^2)^n}{p^{n-1}(1-\lambda p^2)}\quad \mbox{for}\quad n\geq3.
 $$
In the present article, we first obtain a representation formula for functions in the class $\mathcal{V}_p(\lambda)$.
Using this representation, we prove the aforementioned conjecture for $n=3,4,5$ whenever $p$ belongs to certain subintervals of $(0,1)$.
Also we determine non sharp bounds for $|a_n(f)|,\,n\geq 3$ and for $|a_{n+1}(f)-a_n(f)/p|,\,n\geq 2$.
\end{abstract}
\thanks{}
\maketitle
\pagestyle{myheadings}
\markboth{B. Bhowmik and F. Parveen}{On the Taylor coefficients of the class $\mathcal{V}_{p}(\lambda)$}

\bigskip

\section{Introduction}

We shall use following notations throughout the discussion of this article. Let $\IC$ be the whole complex plane, $\ID:=\{z\in \IC: |z|<1\}$ and $\Delta:=\{\zeta\in \IC:|\zeta|>1\}\cup\{\infty\}$. Let $\mathcal{A}$ be the class of all analytic functions $f$ defined in $\ID$ with the normalization $f(0)=0=f'(0)-1$ and $\mathcal{S}=\{f\in \mathcal{A}: f \,\mbox{is \,univalent}\}$. Each $f\in \mathcal{S}$ has the following Taylor expansion:
\beq\label{fp5eq8}
f(z)=z+\sum_{n=2}^{\infty} a_n(f) z^n,\quad z\in \ID.
\eeq
In the last century, the field of geometric function theory provided many interesting and fascinating facts.
One of the main problem of this field was the Bieberbach conjecture which was proposed in the year 1916.
This conjecture states that each $f\in \mathcal{S}$ with the expansion (\ref{fp5eq8}) must satisfy the inequality $|a_n(f)|\leq n$ for all $n\geq 2$.
In the year 1985, L. de Branges \cite{db} proved this conjecture. In order to settle the Bieberbach conjecture prior to the
effort made by de Branges,  many subclasses of $\mathcal{S}$ were introduced that are geometric in nature and the conjecture was being proved for these subclasses.
Some of the special subclasses of $\mathcal{S}$ for which this conjecture was settled were the class of convex functions, starlike functions and close to convex functions.

In \cite{aksen}, Aksent\'{e}v proved a sufficient condition for the functions of the form
$$
F(\zeta)=\zeta+\sum_{n=0}^{\infty} b_{n}\zeta^{-n},\quad \zeta\in\Delta,
$$
to be univalent in $\Delta$. This condition enables many authors to consider the following class of functions:
$$
\mathcal{U}(\lambda):=\left\{ f\in \mathcal{A}: |U_f(z)|<\lambda \,~\mbox{for}~z\in \ID\right\}
$$
where $U_f(z):=(z/f(z))^2 f'(z)-1$ and $0<\lambda \leq 1$. It is well-known that each function in $\mathcal{U}(\lambda)$ is univalent
and $\mathcal{U}(\lambda)\subsetneqq \mathcal{S}$. For a detailed study of the class $\mathcal{U}(\lambda)$ one may go through the articles \cite{fourn, opw, opw1, obrauni} and references therein.
In this note, we consider the meromorphic analogues of the classes $\mathcal{A}$, $\mathcal{S}$ and $\mathcal{U}(\lambda)$.
To this end, let  $\mathcal{A}(p)$ be the class  which is defined as the collection of functions in $\ID$ having a simple pole at $z=p$ where $p\in (0,1)$ and analytic in $\ID\setminus\{p\}$ satisfying the normalization $f(0)=0=f'(0)-1$. Let $\Sigma (p):=\{f\in \mathcal{A}(p) :\, f\, \mbox{is univalent}\}$. In \cite{BF-1}, we established the following sufficient condition for functions in $\mathcal{A}(p)$ to be univalent:

\begin{Thm}\label{TheoA}
Let $f\in \mathcal{A}(p)$. If $|\mathcal{U}_f(z)|\leq((1-p)/(1+p))^2$ for $z\in \ID$,
then $f$ is univalent in $\ID$.
\end{Thm}

Motivated by this sufficient condition stated in the above theorem, in \cite{BF-1}  we considered the following subclass of $\Sigma(p)$.
\bdefe
Let $\mathcal{U}_{p}(\lambda)$ be the family of all functions $f \in \mathcal{A}(p)$ such
that $\left|U_f(z)\right|< \lambda ((1-p)/(1+p))^2$, $z\in \ID$ holds for some $0<\lambda \leq1$.
\edefe
Interested reader may go through the articles \cite{BF-1} and \cite{BF-2} for many other interesting results for functions in the class $\mathcal{U}_p(\lambda)$. We point out here that in \cite{BF-3}, we improve the sufficient condition in Theorem~A for univalence with replacing the number $((1-p)/(1+p))^2$ by $1$ and subsequently the following class of functions $\mathcal{V}_p(\lambda)$ was introduced:
$$
\mathcal{V}_p(\lambda)=\left\{f\in \mathcal{A}(p):\left|U_{f}(z)\right|<\lambda, ~z\in \ID\right\},\quad \mbox{for}~~ \lambda\in(0,1].
$$
In \cite{BF-3}, we proved that $\mathcal{U}_{p}(\lambda)\subsetneq \mathcal{V}_{p}(\lambda)$ and discussed many other aspects of this class of functions. Let $\mathcal{B}$ be the class of functions $w$ which are analytic in $\ID$ and for $z\in \ID$, $|w(z)|\leq1$. In \cite{BF-3}, we proved the following integral representation formula for functions in $\mathcal{V}_p(\lambda)$, i.e., each function in $\mathcal{V}_p(\lambda)$ can be expressed as:
\beq\label{fp5eq6}
\frac{z}{f(z)}=1-\left(\frac{f''(0)}{2}\right)z+\lambda z\int_{0}^{z}w(t)dt,
\eeq
where $w\in\mathcal{B}$. Since each $f\in \mathcal{V}_p(\lambda)$ is analytic in the disc $\ID_p:=\{z:|z|<p\}$,
therefore it has the Taylor expansion of the form (\ref{fp5eq8}) valid in $\ID_p$. In \cite[Theorem~5]{BF-3},
the authors of the present article established the exact region of variability of the second Taylor coefficients $a_{2}(f)$, $f\in \mathcal{V}_p(\lambda)$
which we state below:
$$
|a_{2}(f)-1/p| \leq \lambda p,
$$
and made the following conjecture about the exact bounds for the modulus of the $n$-th Taylor coefficients:

\bcon\label{fp5con2}
If $f\in\mathcal{V}_{p}(\lambda)$ for some $0<\lambda\leq 1$ and has the expansion of the form $(\ref{fp5eq8})$ in $\ID_p$.
Then
\beq\label{fp5eq7}
|a_n(f)|\leq \frac{1-(\lambda p^2)^{n}}{p^{n-1}(1-\lambda p^2)},
\eeq
for $n\geq 3$ and equality occurs in the above inequality for the following functions:
\beq\label{fp5eq15}
k_p^{\lambda}(z):=\frac{-pz}{(z-p)(1-\lambda pz)}.
\eeq
\econ

\br\label{fp5r1}
It is easy to check that as $p\ra 1-$, the inequality (\ref{fp5eq7}) gives the conjectured bound $|a_n(f)|\leq \sum_{k=0}^{n-1}\lambda^k$, $n\geq 3$
for the class $\mathcal{U}(\lambda)$ (see \cite{opw}) and for $\lambda=1$ the above conjecture reduces to the Jenkin's theorem (compare \cite{Jen}) for the class $\Sigma(p)$.
Also by taking $p\rightarrow 1-$ and $\lambda=1$ in  (\ref{fp5eq7}), we will get the famous de Branges theorem for the class $\mathcal{S}$ (compare \cite{db}).
\er

We organise this article as follows. First we prove a representation formula for functions in the class $\mathcal{V}_p(\lambda)$.
Next, with the help of this representation formula we prove the Conjecture~\ref{fp5con2} for $n=3,4,5$ with certain range of values of $p$.
Finally, we obtain non sharp bounds for $|a_n(f)|,\,n\geq 3$ and for $|a_{n+1}(f)-a_n(f)/p|,\,n\geq 2$.

\section{Main Results}
We start this section with the following representation formula for functions in the class $\mathcal{V}_p(\lambda)$:
\bthm
Each $f\in \mathcal{V}_p(\lambda)$ can be represented as
\beq\label{fp5eq2}
f(z)=\frac{-pz}{(z-p)(1-\lambda p z w(z))},\, z\in \ID,
\eeq
where $w\in \mathcal{B}$. Also every $f\in \mathcal{V}_p(\lambda)$ can be expressed as
\beq\label{fp5eq10}
f(z)=\frac{-pz u(z)}{(z-p)(1-\lambda p)},\, z\in \ID,
\eeq
where $u\in \mathcal{B}$ and $u(0)=1-\lambda p$.
\ethm
\bpf
For every $f\in \mathcal{V}_p(\lambda)$ we have from (\ref{fp5eq6}),
$$
\frac{z}{f(z)}=1-\left(\frac{f''(0)}{2}\right)z+\lambda z\int_{0}^{z}w_1(t)dt,\quad z\in \ID,
$$
where $w_1 \in \mathcal{B}$. Since $f(p)=\infty$, the above equality yields
$$
\frac{f''(0)}{2}=\frac{1}{p}\left(1+\lambda p\int_{0}^{p}w_1(t)dt\right)
$$
and hence
\beq\label{fp5eq12}
\frac{z}{f(z)}=1-\frac{z}{p}\left(1+\lambda p\int_{0}^{p}w_1(t)dt\right)+\lambda z\int_{0}^{z}w_1(t)dt,\quad z\in \ID.
\eeq
Let us define
$$
w(z):=\left(\int_{p}^{z}w_1(t)dt\right)/(z-p),\, z\in \ID.
$$
Now it is a simple exercise to see that $|w(z)|\leq 1$ and $w(p)=w_1(p)$. Consequently (\ref{fp5eq12}) takes the following form:
\beqq
\frac{z}{f(z)}=\frac{-(z-p)(1-\lambda p z w(z))}{p},
\eeqq
where $w\in \mathcal{B}$. This proves the representation formula (\ref{fp5eq2}). Next we see that
$$
|1-\lambda p z w(z)|\geq 1-\lambda p |zw(z)|\geq 1-\lambda p.
$$
We now define
\beq\label{fp5eq13}
u(z):=\frac{(1-\lambda p)}{(1-\lambda p z w(z))},\, z\in \ID.
\eeq
Then clearly $|u(z)|\leq 1$ and $u(0)=1-\lambda p$. Now plugging (\ref{fp5eq13}) in (\ref{fp5eq2}) we get (\ref{fp5eq10}). This completes the proof of Theorem~1.
\epf

The next theorem deals with the estimate of $|a_n(f)|$ for $n=3, 4, 5$ under some restriction on the range of values of $p$ where $a_n(f)$ is defined by (\ref{fp5eq8}).

\bthm
Let $f\in \mathcal{V}_p(\lambda)$ have expansion of the form $(\ref{fp5eq8})$.
Then the inequality $(\ref{fp5eq7})$ holds for $n=3,\,p\in (0,1/2]$; for\, $  n=4,\, p\in (0, (\sqrt3-1)/2]$ and for $n=5,\, p\in (0, (\sqrt5-1)/4]$.
Equality holds in the above inequality for the function $(\ref{fp5eq15})$.
\ethm

\bpf
Let each $w\in \mathcal{B}$ has the following Taylor expansion in $\ID$:
$$
w(z)=\sum_{n=0}^{\infty} c_n z^n.
$$
Now inserting the above expression for $w$ and the series expansion $(\ref{fp5eq8})$ for $f$ in the representation formula (\ref{fp5eq2}) we get,
\beqq
z+\sum_{n=2}^{\infty} a_n(f) z^n = z(1-z/p)^{-1}\left(1-\lambda p z \sum_{n=0}^{\infty} c_n z^n\right)^{-1}.
\eeqq
Next comparing the coefficients of $z^n$, $n=3, 4,5$ in the above equality we get
\beq
a_3(f)&=&\lambda p c_1+ \lambda^2 p^2 c_0^2+ \lambda c_0 +1/p^2 \label{fp5eq3},\\
a_4(f)&=&\lambda p c_2+2\lambda^2 p^2c_0 c_1+\lambda^3 p^3 c_0^3+\lambda c_1+\lambda^2 p c_0^2+\lambda c_0/p+ 1/p^3 \label{fp5eq4} \quad \mbox{and}\\
a_5(f)&=&\lambda p c_3+\lambda ^2 p^2 c_1^2+2c_0 c_2 \lambda^2 p^2+3 c_0^2 c_1 \lambda ^3 p^3 +\lambda ^4 p^4 c_0^4+\lambda c_2 \label{fp5eq5}\\
       &&+2c_0 c_1\lambda ^2 p + \lambda ^3 p^2 c_0^3 +\lambda c_1/p+ \lambda^2 c_0 ^2+\lambda c_0/p^2+1/p^4. \nonumber
\eeq
Now from \cite{rus}, we know that
\beq\label{fp5eq9}
|c_0|\leq 1 \quad \mbox{and} \quad |c_n|\leq 1- |c_0|^2 \quad \mbox{for ~all} \quad n\geq 1.
\eeq
Using these inequalities in (\ref{fp5eq3}) we have
\beqq
|a_3(f)|&\leq &\lambda p|c_1|+ \lambda^2 p^2 |c_0|^2+\lambda |c_0|+1/p^2\\
     &\leq & \lambda p (1-|c_0|^2)+\lambda |c_0|+\lambda^2 p^2|c_0|^2+1/p^2.
\eeqq
Setting $|c_0|=x$, we consider $h(x):=\lambda p (1-x^2)+\lambda x+\lambda^2 p^2 x^2+1/p^2$. So $ x \in [0,1]$ and
$$
h'(x)=\lambda (1-2p x)+2 \lambda^2 p^2 x.
$$
Now $h'(x)\geq 0$ for $p\in (0,1/2]$ and for $0<\lambda\leq 1$. This shows that the function $h$ is increasing in $[0,1]$. Therefore,
$$
\max_{0\leq x\leq 1}h(x)=h(1)=\lambda+\lambda^2 p^2+1/p^2
$$
and thus
\beqq
|a_3(f)|\leq 1/p^2+\lambda +\lambda ^2 p^2= (1-(\lambda p^2)^3)/p^2(1-\lambda p^2).
 \eeqq
Using this similar idea, now we prove the conjectured bound for $|a_4(f)|$ and $|a_5(f)|$ when $p$ lies in the intervals stated in the theorem. If we use triangle inequality and the bounds for $|c_n|, n\geq 1$ in (\ref{fp5eq4}), we get
\beqq
|a_4(f)|&\leq& (-2\lambda^2 p^2+\lambda^3 p^3)|c_0|^3+(-\lambda p -\lambda +\lambda ^2 p)|c_0|^2 \\
&&+(2\lambda^2 p^2+\lambda /p)|c_0|+\lambda p+\lambda+1/p^3.
\eeqq
As before setting $|c_0|=x$, we introduce the following function:
\beqq
g(x):=(-2\lambda^2 p^2+\lambda^3 p^3)x^3+(-\lambda p -\lambda +\lambda ^2 p)x^2
+(2\lambda^2 p^2+\lambda /p)x+\lambda p+\lambda+1/p^3.
\eeqq
Since $(-2\lambda^2 p^2+\lambda^3 p^3)<0$ and $(-\lambda p -\lambda +\lambda ^2 p)<0$, then we have
\beqq
g'(x)&=&3(-2\lambda^2 p^2+\lambda^3 p^3)x^2+2(-\lambda p -\lambda +\lambda ^2 p)x+2\lambda^2 p^2+\lambda /p\\
     &\geq& 3(-2\lambda^2 p^2+\lambda^3 p^3)+2(-\lambda p -\lambda +\lambda ^2 p)+2\lambda^2 p^2+\lambda /p \\
     &=& -4 \lambda^2 p^2+3\lambda^3 p^3-2 \lambda p-2 \lambda +2 \lambda ^2 p+ \lambda/p\\
     &=& 3\lambda^3 p^3+ 2\lambda ^2 p(1-2p)+ \lambda(1/p -2-2p).
\eeqq
Now it is a simple exercise to check that $(1-2p)>0$ and $(1/p -2-2p)>0$ for $p\in (0, (\sqrt3-1)/2]$. Hence $g'(x)\geq 0$ for all $0\leq x\leq 1$ and $p\in (0, (\sqrt3-1)/2]$. This implies that $g$ is increasing in $[0,1]$ and
$$
\max_{0\leq x \leq1}g(x)=g(1)=(1+\lambda p^2+\lambda^2 p^4+ \lambda^3 p^6)/p^3.
$$
This settles our claim for the bound of the coefficients $a_4(f)$. Finally we turn our attention into proving the coefficient bound of $|a_5(f)|$. Use of the estimate (\ref{fp5eq9}) in (\ref{fp5eq5}) gives
\beqq
|a_5(f)|&\leq& (\lambda^4 p^4+\lambda^2 p^2-3\lambda^3 p^3) |c_0|^4+(-2 \lambda^2 p^2-2 \lambda^2 p+\lambda^3 p^2)|c_0|^3 \\
          &&+(-\lambda p-2 \lambda^2 p^2+3\lambda^3 p^3-\lambda+\lambda^2-\lambda/p)|c_0|^2\\
          &&+(2\lambda^2 p^2+2\lambda^2 p+\lambda/p^2)|c_0|+\lambda p+\lambda^2 p^2+\lambda+\lambda/p+1/p^4.
\eeqq
Now assume as before $|c_0|=x\in [0,1]$ and define
\beqq
q(x)&:=&(\lambda^4 p^4+\lambda^2 p^2-3\lambda^3 p^3) x^4+(-2 \lambda^2 p^2-2 \lambda^2 p+\lambda^3 p^2)x^3 \\
&&+(-\lambda p-2 \lambda^2 p^2+3\lambda^3 p^3-\lambda+\lambda^2-\lambda/p)x^2 \\
&& +(2\lambda^2 p^2+2\lambda^2 p+\lambda/p^2)x+\lambda p+\lambda^2 p^2+\lambda+\lambda/p+1/p^4.
\eeqq
Thus we have
\beqq
q'(x)&=&4(\lambda^4 p^4+\lambda^2 p^2-3\lambda^3 p^3)x^3+3(-2 \lambda^2 p^2-2 \lambda^2 p+\lambda^3 p^2)x^2 \\
&&+2(-\lambda p-2 \lambda^2 p^2+3\lambda^3 p^3-\lambda+\lambda^2-\lambda/p)x+2\lambda^2 p^2+2\lambda^2 p+\lambda/p^2.
\eeqq
Since $(\lambda^4 p^4+ \lambda^2 p^2)\geq 2(\lambda ^6 p^6)^{1/2}$ and both the quantities $(-2 \lambda^2 p^2-2 \lambda^2 p+\lambda^3 p^2)$ and $(-\lambda p-2 \lambda^2 p^2+3\lambda^3 p^3-\lambda+\lambda^2-\lambda/p)$ are negative, so
\beqq
q'(x)&\geq& 2\lambda^3 p^3 -8\lambda^2 p^2 -4 \lambda^2 p+3 \lambda^3 p^2-2 \lambda p-2 \lambda-2\lambda/p +2 \lambda^2 +\lambda/p^2 \\
     &=&  2\lambda^3 p^3+3 \lambda^3 p^2-2\lambda^2 (4p^2+2p-1)-\lambda(2p +2+2/p-1/p^2).
\eeqq
Again since $(4p^2+2p-1)<0$ and $(2p +2+2/p-1/p^2)<0$ for $p\in (0, (\sqrt5-1)/4]$, the above inequality gives $q'(x)\geq 0$ i.e. , the function $q$ is increasing in the interval $[0,1]$ and hence
$$
\max_{0\leq x \leq 1}q(x)=q(1)=(1+\lambda p^2+\lambda^2 p^4+ \lambda^3 p^6+\lambda^4 p^8)/p^4.
$$
This proves our claim regarding the coefficient bound of $a_5(f)$. Also it is easy to check that
$$
k_p^{\lambda}(z)=\sum_{n=1}^{\infty}\frac{1-\lambda^n p^{2n}}{p^{n-1}(1-\lambda p^2)}z^n, \quad z\in \ID_p.
$$
This proves the sharpness part of the inequality stated in the theorem.
\epf

\br
Using similar lines of proof as in the previous theorem one may be able to prove the conjectured bound for $n\geq 6$ with certain range of values of $p$. However, we expect that this will need much more effort in calculation. We leave this open to the interested reader.
\er

Here we present the following definition which we need for our further discussion.
\begin{defn}
For any two analytic functions $f$ and $g$,  $f$ is said to be subordinate to the function $g$ if the relation $f(z)=g(w(z)),\, z\in \ID$ holds for any $w\in \mathcal{B}$ with $w(0)=0$. This phenomenon is abbreviated as $f\prec g$.
\end{defn}
In the next result we prove non sharp bounds for the absolute value of the Taylor coefficients $|a_{n}(f)|$ for $n\geq 3$ and $|a_{n+1}(f)-a_n(f)/p| $ for $n\geq 2$ whenever $f\in \mathcal{V}_p(\lambda)$.

\bthm
Let $f\in \mathcal{V}_p(\lambda)$ be of the form $(\ref{fp5eq8})$ in  $\ID_p$. Then for $n\geq3$,
\beq\label{fp5eq11}
|a_n(f)|\leq \frac{1}{p^{n-1}}+\left(\sum_{k=1}^{n-1}\lambda^{2k} p^{2k}\right)^{1/2}\left(\sum_{k=1}^{n-1}\frac{1}{p^{2(n-k-1)}}\right)^{1/2},
\eeq
and for $n\geq 2$, we have
$$
|a_{n+1}(f)-a_n(f)/p|\leq \lambda p.
$$
\ethm

\bpf
From the representation (\ref{fp5eq2}) we see that
\beq\label{fp5eq14}
f(z)=f_1(z)f_2(z),\quad z\in \ID,
\eeq
where
$$
f_1(z):=\frac{-pz}{(z-p)} \quad \mbox{and}\quad f_2(z):=\frac{1}{(1-\lambda p z w(z))},\,w\in \mathcal{B}.
$$
Let us now consider the Taylor expansions for the following functions in $\ID$:
$$
f_1(z)=\sum_{n=1}^{\infty}\frac{1}{p^{n-1}}z^n:=\sum_{n=1}^{\infty}A_n z^n,\, f_2(z)=\sum_{n=0}^{\infty}B_n z^n \, \mbox{and} \quad 1/(1-\lambda p z)=\sum_{n=0}^{\infty}(\lambda p)^n z^n.
$$
Since $f_2 \prec 1/(1-\lambda p z)$, therefore using Rogosinski's theorem \cite[Theorem II]{rogo}, we get
\beq\label{fp5eq1}
\sum_{k=1}^{n-1}|B_k|^2\leq \sum_{k=1}^{n-1} \lambda ^{2k} p^{2k}.
\eeq
Now equating coefficients of $z^n$ on both sides of the equation (\ref{fp5eq14}), we get
$$
a_n(f)=A_n+\sum_{k=1}^{n-1}B_k A_{n-k},
$$
where  $A_n=1/p^{n-1}$. Thus an application of the Cauchy-Schwarz inequality yields
\beqq
|a_n(f)|&\leq& |A_n|+\left|\sum_{k=1}^{n-1}B_k A_{n-k}\right| \\
     &\leq& |A_n|+\left(\sum_{k=1}^{n-1}|B_k|^2\right)^{1/2}\left(\sum_{k=1}^{n-1}|A_{n-k}|^2\right)^{1/2} \\
     &\leq& \frac{1}{p^{n-1}}+ \left(\sum_{k=1}^{n-1}\lambda ^{2k} p^{2k}\right)^{1/2}\left(\sum_{k=1}^{n-1}\frac{1}{p^{2(n-k-1)}}\right)^{1/2} \quad \left(\mbox{by}\,(\ref{fp5eq1})\right)
\eeqq
which is the required bound that we wish to prove. Again using the representation (\ref{fp5eq2}) in Theorem~1, we see that
$(z-p)f(z)/(-pz)\prec 1/(1-\lambda p z)=\sum_{n=0}^{\infty}(\lambda p)^n z^n$. Now since $f$ has expansion of the form (\ref{fp5eq8}), we compute
$$
\frac{(z-p)f(z)}{-pz}=1+\sum_{n=1}^{\infty}\left(a_{n+1}(f)-a_n(f)/p\right)z^n.
$$
Next we note that the sequence $\{(\lambda p)^n\}$ is a nonnegative, decreasing and convex sequence. Here we clarify that a real sequence $\{x_n\}$ is called convex sequence if $x_{n-1}+x_{n+1}\geq 2 x_n$ holds for all $n$. Now an application of a well-known result of Rogosinski (see \cite[Theorem VII]{rogo}), we have
$$
|a_{n+1}(f)-a_n(f)/p|\leq \lambda p, \quad \forall\, n\geq  2.
$$
This completes the proof of the theorem.
\epf

\br
In \cite{PW1}, the authors proved the following non sharp bound for $|a_n(f)|,\,n\geq3$ for functions in the class $\mathcal{U}(\lambda)$:
$$
|a_n(f)|\leq 1+\lambda \sqrt{n-1}\sqrt{\sum_{k=0}^{n-2}\lambda ^{2k}}.
$$
We  see that as $ p\rightarrow 1-$, the obtained bound (\ref{fp5eq11}) coincides with the above bound.
\er

\end{document}